\documentclass[11pt]{article}

\usepackage[numbers]{natbib}
\usepackage{amsmath}
\usepackage{amssymb}
\usepackage{amsthm}
\usepackage{xcolor}
\usepackage{graphicx}
\usepackage{tikz}
\usetikzlibrary{arrows.meta}
\usepackage{float}
\usepackage{array}
\usepackage{booktabs}
\usepackage{tabularx}
\usepackage{hyperref}
\usepackage{algorithm}
\usepackage{algpseudocode}
\usepackage{enumitem}\usepackage[margin=1in]{geometry}

\hypersetup{colorlinks=true,citecolor=blue,linkcolor=blue,hypertexnames=false}
\allowdisplaybreaks

\newcommand{\paperTitle}{A Linear Bound for the Affine Banach–Mazur Diameter}
\newcommand{\paperAuthor}{Zhao Song\thanks{\texttt{magic.linuxkde@gmail.com}.} \\ Independent Researcher \and Xinzhi Zhang\thanks{\texttt{xinzhizhang@microsoft.com}.}\\ Microsoft Research}

\newcommand{\E}{\mathbb E}
\newcommand{\R}{\mathbb R}

\renewcommand{\d}{\mathrm{d}}

\theoremstyle{plain}
\newtheorem{theorem}{Theorem}[section]
\newtheorem{lemma}[theorem]{Lemma}
\newtheorem{fact}[theorem]{Fact}

\theoremstyle{definition}
\newtheorem{definition}[theorem]{Definition}

\begin{document}

\date{November 27, 2026}
\title{\paperTitle}
\author{\paperAuthor}
\maketitle

\begin{abstract}
The affine Banach--Mazur diameter $D_n$ is the supremum of the affine Banach--Mazur distances between convex bodies in $\mathbb{R}^n$. In this paper, we resolve the problem of determining the asymptotic order of $D_n$ by proving that $D_n\le 9.6 n$ for every $n\ge 2$. This matches the lower bound $D_n\ge n$~\cite{lei59} with up to a constant factor.

\end{abstract}

\newpage

\section{Introduction}

Let $n\ge2$, for convex bodies $K,L\subset\mathbb R^n$, define the affine Banach–Mazur distance $d_{\mathrm{BM}}(K,L)$ between $K$ and $L$ as the infimum of all $\lambda\ge1$ for which there exist invertible affine maps $F,G:\mathbb R^n\to\mathbb R^n$ such that
\[
    F(K)\subseteq G(L)\subseteq\lambda F(K).
\]
To quantify how far apart convex bodies can be in this distance, we consider the affine Banach–Mazur diameter in dimension \(n\),
\[
    D_n:=\sup_{K,L}d_{\mathrm{BM}}(K,L),
\]
where the supremum is taken over all convex bodies \(K,L\subset\mathbb R^n\). A natural problem is to determine how \(D_n\) grows with the dimension \(n\) as $n\to \infty$.
 
The first upper bound on $D_n$ follows from John's ellipsoid theorem~\cite{joh48}, proved in 1948: every $n$-dimensional convex body contains an ellipsoid whose dilation by a factor of $n$ about its center contains the body.
Every body is therefore within distance $n$ of the Euclidean ball, and the multiplicative triangle inequality yields $D_n\le n^2$. In the other direction, Leichtweiss~\cite{lei59} showed in 1959 that the $n$-simplex is at distance exactly $n$ from the ball, so $D_n\ge n$. In 1998, Litvak and Tomczak-Jaegermann improved the upper bound to $O(n^{3/2}\sqrt{\log n})$, as in~\cite{blps99}. Rudelson~\cite{rud00} improved this bound to $D_n=O(n^{4/3}\log^9 n)$. In 2025, Bizeul and Klartag proved $D_n=O(n\log^4 n)$ in the first version of their paper~\cite{bk25}. In the third version (2026), they improved the bound to $D_n=O(n\log^{7/2}n)$~\cite{bk25}. Bizeul~\cite{biz26} further improved the bound to $D_n=O(n\log n)$.

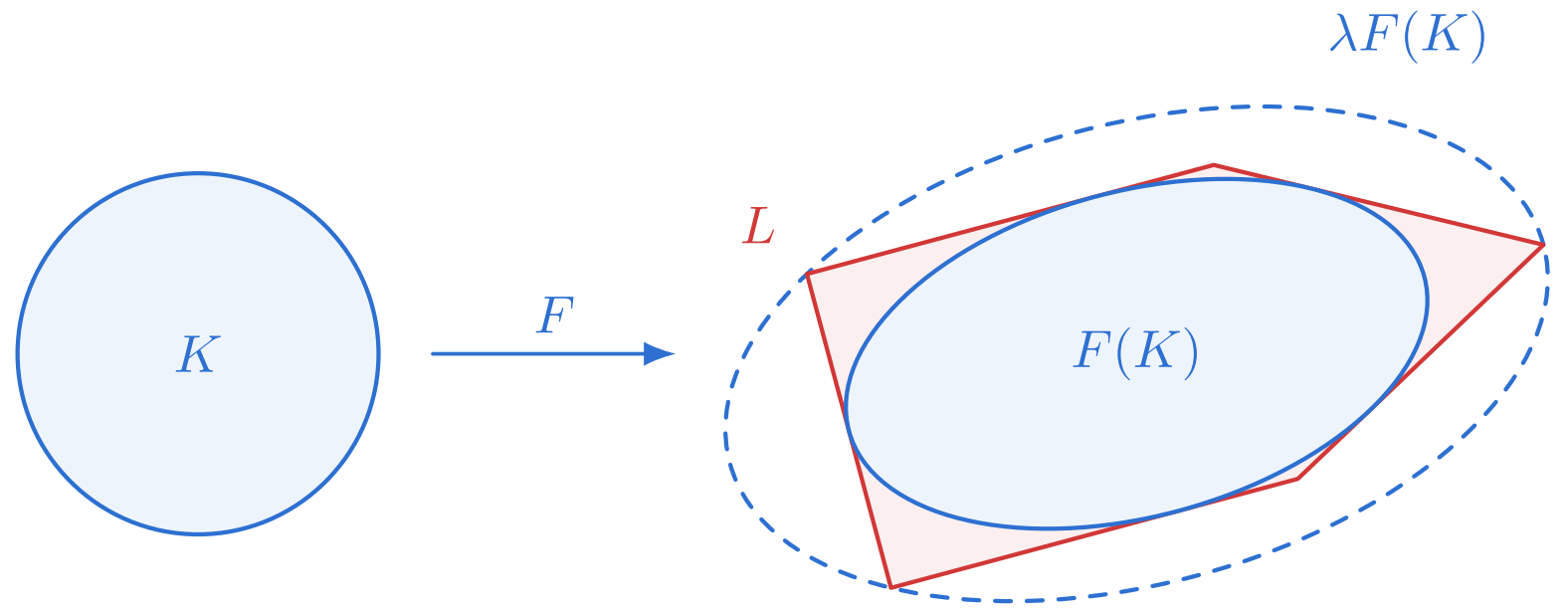
\begin{figure}[H]
    \centering
    \definecolor{bmblue}{RGB}{45,112,209}
    \definecolor{bmred}{RGB}{210,55,55}
    \begin{tikzpicture}[scale=1.2,line cap=round,line join=round,>=Latex,font=\small]
        \draw[bmblue,fill=bmblue!8,line width=0.8pt]
            (-3.4,0) circle[radius=1cm];
        \node[text=bmblue] at (-3.4,0) {$K$};
        \draw[->,bmblue,line width=0.7pt] (-2.1,0) -- node[above,text=bmblue] {$F$} (-0.75,0);
        \begin{scope}[shift={(1.8,0)},rotate=15]
            \pgfmathsetmacro{\bmoutera}{1.65*sqrt(2)}
            \pgfmathsetmacro{\bmouterb}{0.9*sqrt(2)}
            \draw[bmblue,dashed,line width=0.8pt]
                (0,0) ellipse[x radius=\bmoutera cm,y radius=\bmouterb cm];
            \draw[bmred,fill=bmred!8,line width=0.8pt]
                (-1.65,-0.9) -- ({1.65*(sqrt(2)-1)},-0.9)
                -- ({1.65*sqrt(2)},0)
                -- ({1.65*(sqrt(2)-1)},0.9) -- (-1.65,0.9) -- cycle;
            \draw[bmblue,fill=bmblue!8,line width=0.8pt]
                (0,0) ellipse[x radius=1.65cm,y radius=0.9cm];
            \node[text=bmblue] at (0,0) {$F(K)$};
            \node[above left=2pt,text=bmred] at (-1.65,0.9) {$L$};
            \node[anchor=south west,text=bmblue] at (1.3,1.2) {$\lambda F(K)$};
        \end{scope}
    \end{tikzpicture}
    \caption{An illustration of the affine Banach--Mazur distance.}
    \label{fig:banach_mazur_distance}
\end{figure}

\subsection{Our results}

Our main result is stated as follows:
\begin{theorem}[Main result, informal version of Theorem \ref{thm:main}]\label{thm:main_informal}
For every $n\ge2$, the affine Banach--Mazur diameter of $n$-dimensional convex bodies satisfies $D_n\le9.6n$.
\end{theorem}
The lower bound for this problem given is $D_n \geq n$, as shown in \cite{lei59}. We believe that finding the right constant for this problem is an interesting further direction.

\begin{table}[!ht]
    \centering
    \begin{tabular}{|l|l|l|l|}
        \hline
        Year & Author & Reference & Upper/lower bound \\
        \hline
        1948 & John & \cite{joh48} & $O(n^2)$ \\
        1959 & Leichtweiss & \cite{lei59} & $\Omega(n)$ \\
        1998 & Litvak and Tomczak-Jaegermann & recorded in~\cite{blps99} & $O(n^{3/2}\sqrt{\log n})$ \\
        1998 & Rudelson & \cite{rud00} & $O(n^{4/3}\log^9 n)$ \\
        2025 & Bizeul and Klartag & \cite[v1]{bk25} & $O(n\log^4 n)$ \\
        2026 & Bizeul and Klartag & \cite[v3]{bk25} & $O(n\log^{7/2}n)$ \\
        2026 & Bizeul & \cite{biz26} & $O(n\log n)$ \\
        2026 & This paper & Theorem~\ref{thm:main_informal} & $O(n)$ \\
        \hline
    \end{tabular}
    \caption{Overview of the previous Banach--Mazur diameter bounds for arbitrary convex bodies.}
    \label{tab:banach_mazur_history}
\end{table}

\section{Preliminaries}
\label{sec:affine_preliminaries}

For positive integer $n$, we define set $[n] := \{1,2,\cdots,n\}$. For a vector $x$, we use $x^\top$ to denote its transpose.

\begin{definition}[Affine Banach--Mazur distance and diameters; see~\cite{las07}]\label{def:affine_distance}
A convex body is a compact convex set with nonempty interior. For two
convex bodies $K,L\subset\R^n$, write
\[
 d_{\mathrm{BM}}(K,L):=\inf\{\lambda\ge1:B\subseteq U\subseteq\lambda B,
 \ B,U\text{ are invertible affine images of }K,L\text{ respectively}\}.
\]
For $n\ge2$, let $K,L$ range over convex bodies in $\R^n$ and $S,T$ over
centrally symmetric convex bodies in $\R^n$. Define
\[
 D_n:=\sup_{K,L}d_{\mathrm{BM}}(K,L),\quad  D_n^{\mathrm{mix}}:=\sup_{K,S}d_{\mathrm{BM}}(K,S),\quad  D_n^{\mathrm{sym}}:=\sup_{S,T}d_{\mathrm{BM}}(S,T).
\]
\end{definition}

We first check that either body may serve as the inner affine image when estimating the distance.

\begin{fact}[Symmetry of the distance]\label{fac:affine_distance_symmetry}
For convex bodies $K,L\subset\R^n$, $d_{\mathrm{BM}}(K,L)=d_{\mathrm{BM}}(L,K)$.
\end{fact}
\begin{proof}
For every admissible pair we have $\lambda^{-1}U\subseteq B\subseteq U$,
where the first inclusion follows by scaling $U\subseteq\lambda B$ by
$\lambda^{-1}$, and the second is the original inner inclusion. Thus the
same $\lambda$ is admissible with the bodies interchanged. Taking infima and interchanging $K,L$ gives the claim.
\end{proof}

We describe set containments through linear inequalities, using support functions to normalize their right-hand sides.

\begin{definition}[Support function]\label{def:affine_support_function}
For a convex body $K$, its support function is
$h_K(b):=\max\{b^\top x : x\in K\}$.
\end{definition}

The following characterization will convert bounds in every normalized support direction into a containment.

\begin{lemma}[Normalized support half-spaces]\label{lem:affine_support}
If $0\in\operatorname{int} K$, then
$K=\bigcap_{h_K(b)=1}\{x:b^\top x\le1\}$.
\end{lemma}
\begin{proof}
Every point of $K$ satisfies these inequalities. For $x\notin K$, choose
a nearest point $y\in K$ and set $u:=x-y$. Minimizing squared distance
along each segment from $y$ to $p\in K$ gives $u^\top(p-y)\le0$.
Thus $h_K(u)=u^\top y<u^\top x$, where the first step follows because
$y$ attains the maximum of $u^\top p$ over $p\in K$, and the second
follows from $u^\top(x-y)=\|x-y\|^2>0$. Since $K$ contains a ball about $0$,
$h_K(u)>0$, and $b:=u/h_K(u)$ excludes $x$.
\end{proof}

We will use the following compactness fact to prove that the objective in Definition~\ref{def:affine_objective} attains its maximum.

\begin{lemma}[Finite convex combinations]\label{lem:affine_finite_hull}
If $S\subset\R^m$ is nonempty and compact, its convex hull is compact
and each of its points is a convex combination of at most $m+1$ members of $S$.
\end{lemma}
\begin{proof}
Given a convex combination with more than $m+1$ positive weights $p_i$,
the vectors $(1,s_i)$ are linearly dependent. Take a nonzero dependence
$\sum_i u_i(1,s_i)=0$ with some $u_i>0$. Set
$t:=\min_{u_i>0} p_i/u_i$. Replacing $p_i$ by $p_i-tu_i$
preserves the point and the total weight, and eliminates at least one term.
Iterate. The hull is therefore the continuous image of $S^{m+1}$ times
the compact simplex of weights under the convex-combination map.
\end{proof}

Sylvester~\cite{syl83} observed that for square matrices $A$ and $B$
of the same order, $AB$ and $BA$ have the same eigenvalues.
The following rectangular version follows by evaluating the
characteristic-polynomial identity of Horn and
Johnson~\cite[Theorem~1.3.22, p.~65]{hj13} at $-1$.

\begin{lemma}[Sylvester's determinant identity]\label{lem:affine_sylvester}
For $U,V\in\R^{n\times k}$, we have
\[
  \det(I_n+UV^\top)=\det(I_k+V^\top U).
\]
\end{lemma}
\begin{proof}
In $(\begin{smallmatrix} I_n&U\\-V^\top&I_k\end{smallmatrix})$,
adding $V^\top$ times the top block row to the bottom block row gives
determinant $\det(I_k+V^\top U)$. Subtracting $U$ times the original
bottom block row from the original top block row gives $\det(I_n+UV^\top)$.
\end{proof}

We use Pr\'ekopa's theorem to prove that the tail of an affine function under a log-concave probability measure is log-concave. This property is used to bound exponential moments in Lemma~\ref{lem:affine_exponential_moment}.

\begin{lemma}[Pr\'ekopa's marginal theorem; Theorem~6 of~\cite{pre73}]
\label{lem:affine_prekopa}
Let $f:\R^d\times\R^m\to[0,\infty)$ be bounded, measurable and
log-concave. Suppose $f(x,y)$ vanishes outside a fixed compact set of
$x$-values. Then $g(y):=\int_{\R^d} f(x,y)\,\d x$ is log-concave.
\end{lemma}

We provide the proof of this lemma for completeness.

\begin{proof}
It suffices to integrate one coordinate at a time. For the one-dimensional
step fix $y_0,y_1$ and $0<\alpha<1$. Set $f_i(u):=f(u,y_i)$,
$I_i:=\int f_i(u)\,\d u$ for $i=0,1$, and
$y_\alpha:=(1-\alpha)y_0+\alpha y_1$.
If either mass is zero the desired inequality holds. Otherwise each $f_i$
is positive and continuous on the interior of its interval of positivity,
by concavity of $\log f_i$ there. 
Define
$F_i(u):=I_i^{-1}\int_{-\infty}^{u} f_i(v)\,\d v$ as the normalized
distribution functions. On those interiors they are strictly increasing
onto $(0,1)$. Set $T(u):=F_1^{-1}(F_0(u))$. Differentiating
$F_1(T(u))=F_0(u)$ and using $\frac{\d F_i}{\d u}(u)=f_i(u)/I_i$ gives
$\tfrac{\d T}{\d u}(u)=I_1 f_0(u)/(I_0 f_1(T(u)))$.
For $w(u):=(1-\alpha)u+\alpha T(u)$ we have
$\tfrac{\d w}{\d u}(u)=(1-\alpha)+\alpha \tfrac{\d T}{\d u}(u)\ge (\tfrac{\d T}{\d u}(u))^\alpha$, where the first step
is differentiation and the second is the weighted arithmetic--geometric
mean inequality. Let $I$ denote the interior of the positivity interval of $f_0$. Then
\[
\begin{aligned}
 \int_{\R} f(v,y_\alpha)\,\d v
 &\ge \int_I f(w(u),y_\alpha)\tfrac{\d w}{\d u}(u)\,\d u\\
 &\ge \int_I f_0(u)^{1-\alpha} f_1(T(u))^\alpha (\tfrac{\d T}{\d u}(u))^\alpha\,\d u\\
 &= I_0^{1-\alpha}I_1^\alpha,
\end{aligned}
\]
where the first step restricts the nonnegative integral to $w(I)$ and
changes variables by $v=w(u)$, the second applies joint log-concavity at
$(u,y_0)$ and $(T(u),y_1)$ together with the bound on $\d w/\d u$, and the third
substitutes the formula for $\d T/\d u$. Indeed, the resulting integrand is
$(I_1/I_0)^\alpha f_0(u)$, whose integral uses $\int_I f_0(u)\,\d u=I_0$.
Each integration preserves boundedness and the common bounded support
in the remaining integration variables. Iteration proves the claim.
\end{proof}

\section{An upper bound for the affine Banach--Mazur diameter}
\label{sec:affine_upper_bound}

In this section, we prove the following main result.

\begin{theorem}[Main result, formal version of Theorem~\ref{thm:main_informal}]\label{thm:main}
Let $K,L\subset\R^n$ be compact convex sets with nonempty interiors. Then
\[
                 d_{\mathrm{BM}}(K,L)\le \frac{48}{5}n .
\]
Consequently $D_n\le 9.6n$ for unrestricted pairs of convex bodies.
\end{theorem}

\subsection{A maximizing affine map}

We construct affine maps with an explicit forward containment, then control the dilation factor by maximizing a determinant objective.

\begin{definition}[Admissible lists and anchored maps]\label{def:affine_list}
Let $K,L\subset\R^n$ be convex bodies with $0\in\operatorname{int}K$. Write
$\mu(B):=\operatorname{vol}(B\cap K)/\operatorname{vol}(K)$ for the uniform probability of a Borel set $B$ under the law on $K$.  An admissible list of length $k$ consists of pairs $(y_i,d_i)$ such that, for every $i\in[k]$,
the label $y_i$ belongs to $L$ and the affine function $d_i(x):=t_i-b_i^\top x$ is nonnegative on $K$. Set
\[
 D(x):=\sum_i d_i(x),\qquad m:=\sum_i b_i,\qquad H:=\max\{D(x) : x\in K\}.
\]
For an anchor $w\in L$, define $F_w(x):=\sum_i d_i(x)(y_i-w)$ and denote its
linear part by $A_w:=-\sum_i(y_i-w)b_i^\top$.
\end{definition}

The parameter $H$ controls the forward containment of $F_w(K)$ in a dilate of $L-w$.

\begin{lemma}[Forward containment]\label{lem:affine_forward_placement}
For every admissible list in Definition~\ref{def:affine_list} and every $w\in L$,
\begin{equation}\label{eq:affine_forward}
                      F_w(K)\subseteq H(L-w).
\end{equation}
\end{lemma}
\begin{proof}
If $H>0$, then $\sum_i d_i(x)/H=D(x)/H\le1$, where the first step is the
definition of $D$ and the second uses $D(x)\le H$. Assigning the remaining
mass to $0\in L-w$ expresses $F_w(x)/H$ as a convex combination in $L-w$.
If $H=0$, every $d_i$ vanishes on $K$ and the containment also holds.
\end{proof}

We will choose the points $y_i$ and affine functions $d_i$ in Definition~\ref{def:affine_list} by maximizing an objective over all admissible lists. We first introduce a continuous piecewise linear function $J:[0,1]\to[1,\infty)$ used to define the determinant term of this objective. We require $J(0)=1$, $J$ to be nondecreasing, and $t\mapsto J(t)/t$ to be nonincreasing on $(0,1]$. The particular knots, given in Section~\ref{sec:affine_certificate_data}, are also chosen for the uniform scalar bound in Lemma~\ref{lem:affine_scalar_certificate}.

\begin{definition}[Determinant function and scalar bounds]\label{def:affine_parameters}
Set $\gamma:=29/44$, $C:=1392/365$, and $M(s):=(1-s)^{-\gamma}$ for $0\le s<1$.
Let $J:[0,1]\to\R$ be the piecewise linear interpolant of the knots in
Table~\ref{tab:affine_knots}. Write $t_+:=\max\{t,0\}$ for the positive part of a real number $t$.
For $u,v\ge0$ set $B:=\max\{u,v\}$ and $b:=\min(u,v)$. Define
\[
 N(u,v):=B J(b/B)\quad(B>0),\qquad N(0,0):=0.
\]
\end{definition}

We use monotonicity of $N$ to compare bounds on its two arguments separately. We use positive homogeneity to compute how the numerator of the objective changes when all $d_i$ are multiplied by the same factor.

\begin{lemma}[Monotonicity of the determinant function]\label{lem:affine_monotone}
The function $N$ in Definition~\ref{def:affine_parameters} is continuous, symmetric,
positively homogeneous, and nondecreasing in each coordinate.
Moreover, $B\le N(u,v)\le J(1) B$.
\end{lemma}
\begin{proof}
For each pair of adjacent interpolation points \((x_i,j_i)\) and \((x_{i+1},j_{i+1})\) defining \(J\), Table~\ref{tab:affine_knots} gives \(j_i\le j_{i+1}\) and \(x_i j_{i+1}\le x_{i+1}j_i\).
Thus each linear piece has nonnegative slope and nonnegative intercept.
It follows that $J(t)$ is nondecreasing and $J(t)/t$ is nonincreasing for $t>0$.
For fixed $B$, increasing $b$ increases $J(b/B)$ weakly. For fixed $b>0$, increasing
$B$ decreases $b/B$, so $b J(b/B)/(b/B)$ increases weakly. When $b=0$,
$N=B$ since $J(0)=1$. These facts give coordinatewise monotonicity, including when
the two coordinates exchange order. Symmetry and positive homogeneity follow from
the formula. Since $1\le J\le J(1)$, the stated bounds also prove continuity at
the origin. Continuity elsewhere follows from continuity of $J$.
\end{proof}

We now use $N$ to define the objective. We also introduce the probability measure with density proportional to $e^D$ relative to $\mu$, so that changes in $\int_K e^D\,\d\mu$ can be expressed as expectations.

\begin{definition}[Determinant objective and tilted measure]\label{def:affine_objective}
For an admissible list in Definition~\ref{def:affine_list}, let
\[
 U_+:=\max\{(\det A_w)_+ : w\in L\},\qquad U_-:=\max\{(-\det A_w)_+ : w\in L\} .
\]
Using the function $N$ and the constant $\gamma$ from Definition~\ref{def:affine_parameters}, define
\[
 \mathcal N:=N(U_+,U_-),\qquad \Phi:=\mathcal N/(\int_K e^{D(x)} \d\mu(x))^\gamma .
\]
For this list, let $X$ have law $\nu$, where
$\d\nu:=e^D \d\mu/\int_K e^D \d\mu$, and set $\bar x:=\E[X]$.
\end{definition}

We next prove that $\Phi$ attains its maximum at an admissible list with $\mathcal N>0$ and $H=O(n)$. These properties ensure that some $F_w$ is invertible and that its containment in Lemma~\ref{lem:affine_forward_placement} has a dilation factor of order $n$.

\begin{lemma}\label{lem:affine_maximizer}
There is an admissible list maximizing the objective $\Phi$ of
Definition~\ref{def:affine_objective} with $\mathcal N>0$. At any
such maximizer, $44n/29\le H\le 73n/29$.
\end{lemma}
\begin{proof}
Fix $w\in L$ and choose $n$ labels so that the $y_i-w$ form a basis.
Choose a basis $b_1,\ldots,b_n$ and take $t_i:=h_K(b_i)$. The matrix $A_w$
is then invertible. For any admissible list and every $w\in L$,
\[
 |\det A_w|\operatorname{vol}(K)=\operatorname{vol}(F_w(K))
 \le\operatorname{vol}(H(L-w))=H^n\operatorname{vol}(L),
\]
where the first step is the volume formula for an affine image, the second
uses Lemma~\ref{lem:affine_forward_placement}, and the third uses translation
invariance and the scaling of $n$-dimensional volume. If $D(x_*)=H$, then $D\ge H/2$ on $(x_*+K)/2$, a set
of $\mu$-mass $2^{-n}$. Hence $\int_K e^D\d\mu\ge 2^{-n}e^{H/2}$.
The numerator of $\Phi$ grows at most as a constant times $H^n$ by Lemma~\ref{lem:affine_monotone} and the determinant bound, whereas its denominator grows at least as
$(2^{-n}e^{H/2})^\gamma$. Thus $\Phi\to0$ as $H\to\infty$.

Let $\Phi_*>0$ be the value of the admissible list constructed at the start of the proof. Choose $H_0$ so that every list with $H>H_0$ has $\Phi<\Phi_*$, which is possible because $\Phi\to0$ as $H\to\infty$. For any list, omit identically zero functions and set $\tau:=\sum_i t_i$. Since $0\in\operatorname{int}K$, every remaining function has $t_i>0$, and $u_i:=b_i/t_i$ belongs to $K^\circ:=\{u:h_K(u)\le1\}$. Also $\tau=D(0)\le H$, where the equality uses $d_i(0)=t_i$ and the inequality uses the definition of $H$. Thus every list with $\Phi\ge\Phi_*$ has $\tau\le H_0$.

For $\tau\le H_0$, the possible moment vectors form the compact set
\[
 \{\tau v:0\le\tau\le H_0,\quad
 v\in\operatorname{conv}\{(1,u,y,yu^\top):u\in K^\circ,\ y\in L\}\}.
\]
Indeed, the weights $t_i/\tau$ of each nonzero list sum to one. Conversely, by Lemma~\ref{lem:affine_finite_hull}, the convex hull is compact and each $v$ has a finite convex representation. Multiplying its weights by $\tau$ gives an admissible list, since $1-u^\top x\ge0$ on $K$. At $\tau=0$ the zero list suffices. Scalar multiplication is continuous, so the displayed set is compact. The moment coordinates give $m=\sum_i t_i u_i$ and $Z:=\sum_i t_i y_i u_i^\top$. Substitution into Definition~\ref{def:affine_list} gives $D(x)=\tau-m^\top x$ and $A_w=wm^\top-Z$. The denominator of $\Phi$ is positive and continuous in these moments, and the numerator is continuous by continuity of $N$ and maximization over the compact anchor set $L$. Hence $\Phi$ attains a maximum on the displayed set. Its value is at least $\Phi_*$, while every list with $\tau>H_0$ has $H\ge\tau>H_0$ and hence $\Phi<\Phi_*$. This proves global maximality over all finite admissible lists.

Fix a maximizer. Scaling all $d_i$ by $t>0$
multiplies each determinant by $t^n$, and hence multiplies $\mathcal N$ by $t^n$ by positive homogeneity. Let $\Phi_t$ denote the objective of
the scaled list. Differentiating its logarithm gives
\[
 \frac{\d}{\d t}\log\Phi_t|_{t=1}
 =n-\gamma\frac{\int_K D(x)e^{D(x)} \d\mu(x)}
                       {\int_K e^{D(x)} \d\mu(x)}
 =n-\gamma\E[D(X)],
\]
where the first step differentiates $n\log t$ and the logarithm of the
normalizing integral (the function $D$ is bounded on $K$), and the second
uses the tilted measure of Definition~\ref{def:affine_objective}. Hence $\E[D(X)]=44n/29\le H$, where the equality
sets this derivative to zero at the maximum, and the inequality uses
$D(X)\le H$. For $0<v\le1$ the sets
$x_*+v(K-x_*)$ increase with $v$. Writing $K_v:=x_*+v(K-x_*)$, we have
\[
 \int_{K_v}e^{D(y)} \d y
 =v^n\int_K e^{D(x_*+v(x-x_*))} \d x
 =v^n\int_K e^{H+v(D(x)-H)} \d x,
\]
where the first step is the change of variables $y=x_*+v(x-x_*)$ with
Jacobian $v^n$, and the second uses the affinity of $D$ and $D(x_*)=H$.
Set $I(v):=\int_{K_v}e^{D(y)} \d y$. Since the sets increase with $v$,
\[
 0\le \frac{1}{I(1)}\frac{\d I}{\d v}|_{v=1^-}
 = n+\E[D(X)]-H,
\]
where the first step follows from the monotonicity of $I$ and the positivity of $I(1)$,
and the second step follows from applying the product rule to the integral representation of $I(v)$ over $K$
at $v=1$ from the left and using $\d\nu(x)=e^{D(x)}\d x/I(1)$. Since $D$ is bounded on $K$,
we may differentiate under the integral. Thus $H\le n+\E[D(X)]=73n/29$,
where the first step rearranges the derivative inequality and the second
substitutes $\E[D(X)]=44n/29$.
\end{proof}

\subsection{Determinant perturbations and scalar inequalities}
\label{sec:affine_perturbation}

By Lemma~\ref{lem:affine_maximizer}, we can fix an admissible list maximizing $\Phi$, with $\mathcal N>0$ and $H=O(n)$. We now derive scalar inequalities needed for the reverse containment. These inequalities follow by appending a nonnegative affine function and a label in $L$ and comparing the resulting objective with its maximum.

\begin{definition}[Anchor coordinates]\label{def:affine_anchors}
For a maximizing list from Lemma~\ref{lem:affine_maximizer}, choose
$z\in L$ maximizing $|\det A_z|$ and set $A:=A_z$. Some anchor has a
nonzero determinant since $\mathcal N>0$, so $A$ is invertible. Set
\[
 G:=A^{-1}(L-z),\qquad \delta_p:=-m^\top p\quad(p\in G).
\]
Choose $r\in G$ maximizing $\delta_r$. Set $a:=\delta_r-1$, and for
each $p\in G$ set $q:=\delta_p-1$.
\end{definition}

Before perturbing the list, we express the determinant numerator $\mathcal N$ using the anchor parameter $a$ and $|\det A|$.

\begin{lemma}[Determinants at the anchors]\label{lem:affine_anchor_envelope}
In Definition~\ref{def:affine_anchors}, the relative anchor determinants
range over $[-a,1]$, and $-1\le q\le a\le1$. With $J$ from
Definition~\ref{def:affine_parameters},
\begin{equation}\label{eq:affine_weight_envelope}
 \mathcal N=|\det A|J(a_+).
\end{equation}
\end{lemma}
\begin{proof}
Changing the anchor from $z$ to $z+Ap$ adds $Apm^\top$ to $A$.
By Lemma~\ref{lem:affine_sylvester} with $k=1$, the relative determinant is
$\det(I+pm^\top)=1+m^\top p=1-\delta_p$, where the first step reduces the
determinant to a scalar and the second uses the definition of $\delta_p$.
The choice of $z$ gives $|1-\delta_p|\le1$. Thus $-1\le q\le a\le1$,
where the first inequality uses $|1-\delta_p|\le1$, the second the choice
of $r$, and the third $|1-\delta_r|\le1$.
Since $0,r\in G$ and $G$ is convex, $1-\delta_p$ takes every value between
$1$ and $-a$.
The positive and negative extrema of the relative determinants are therefore
$1$ and $a_+$. Multiplying by $\det A$ scales their magnitudes by $|\det A|$
and exchanges their roles when $\det A<0$. Symmetry and homogeneity in
Lemma~\ref{lem:affine_monotone} give
$\mathcal N=|\det A|N(1,a_+)=|\det A|J(a_+)$,
where the last step uses $0\le a_+\le1$ in Definition~\ref{def:affine_parameters}.

\end{proof}

We now construct the nonnegative affine functions used in the perturbation. For a support direction $b$, we use the deficit $h_K(b)-b^\top x$ and choose the normalization so that its expectation under $\nu$ equals one.

\begin{definition}[Normalized support deficit and coordinates]\label{def:affine_deficit}
Use $\bar x$ from Definition~\ref{def:affine_objective} and the anchor
coordinates of Definition~\ref{def:affine_anchors}. For a support direction
$b$ normalized by $h_{K-\bar x}(b)=1$ and a point $p\in G$, set
\[
 T(x):=h_K(b)-b^\top x,\qquad P:=b^\top p,\qquad R:=b^\top r .
\]
Here $h_K$ is the support function in Definition~\ref{def:affine_support_function}.
\end{definition}

We estimate $(\E[e^{sT(X)}])^\gamma$, the factor by which $D\mapsto D+sT$ changes the denominator of $\Phi$, uniformly over normalized support directions.

\begin{lemma}[Exponential moment of the support deficit]\label{lem:affine_exponential_moment}
The barycenter $\bar x$ belongs to $\operatorname{int}K$. For each normalized
$b$ in Definition~\ref{def:affine_deficit}, $T\ge0$ on $K$, $\E[T(X)]=1$,
and $(\E[e^{sT(X)}])^{\gamma}\le M(s)$ for $0<s<1$.
\end{lemma}
\begin{proof}
A closed ball inside $\operatorname{int}K$
carries a positive uniform component of $\nu$, since its density is
bounded below there. Thus $\bar x$ is a convex combination with positive
weight on the center of that ball and remaining weight on a point of $K$.
In particular $\bar x\in\operatorname{int}K$. By the definition of $h_K$, $T\ge0$ and
\[
 \E[T(X)]=h_K(b)-b^\top\bar x=h_{K-\bar x}(b)=1,
\]
where the first step uses linearity of expectation and $\bar x=\E[X]$,
the second translates the support function by $-\bar x$, and the third
uses the normalization of $b$.

Set $S(t):=\Pr[T(X)\ge t]$. The function
$(x,t)\mapsto e^{D(x)}\mathbf{1}_{\{x\in K,\ t\le T(x)\}}$ is log-concave because its support is convex and $D$ is affine there. It is bounded and has
$x$-support in $K$. Integrating in $x$ gives a positive constant times
$S(t)$, so $S$ is log-concave by Lemma~\ref{lem:affine_prekopa}.
Since $S(0)=1$, interpolation between $0$ and $u+v$ at $u$ and at $v$
gives $S(u+v)\le S(u)S(v)$ for $u,v\ge0$. For each integer $k\ge1$,
integration over the nonnegative orthant now yields
\[
 \frac{\E[T(X)^k]}{k!}
 =\int_{[0,\infty)^k}S(t_1+\cdots+t_k) \d t_1\cdots \d t_k
 \le (\int_0^\infty S(t) \d t)^k=1,
\]
where the first step integrates the volume $T(X)^k/k!$ of the simplex
$\{t_i\ge0:\sum_i t_i\le T(X)\}$ and exchanges expectation and integration
by nonnegativity, the second uses $S(\sum_i t_i)\le\prod_i S(t_i)$
and factors the integral, and the third uses the tail integral identity
$\int_0^\infty S(t) \d t=\E[T(X)]=1$. It follows that
\[
 \E[e^{sT(X)}]=\sum_{k\ge0}\frac{s^k\E[T(X)^k]}{k!}
 \le\sum_{k\ge0}s^k=(1-s)^{-1},
\]
where the first step uses monotone convergence for the exponential series,
the second applies the moment bound term by term, and the third sums
the geometric series for $0<s<1$. Raising both sides to the power $\gamma$
gives the claim.
\end{proof}

By Lemma~\ref{lem:affine_exponential_moment}, appending $sT$ with a label in $L$ changes the denominator of $\Phi$ by a factor of at most $M(s)$. Maximality of $\Phi$ implies that the perturbed determinant numerator is at most $M(s)$ times the original numerator.

\begin{definition}[Scalar inequalities]\label{def:affine_rows}
For $-1\le q\le a\le1$ and $P,R\in\R$, require for every $0<s<1$ that
\begin{equation}\label{eq:affine_single_row}
 -q+sP\le J(a_+)M(s).
\end{equation}
For $0\le a\le1$ also require
\begin{equation}\label{eq:affine_rows}
\begin{aligned}
 N(\max\{1,-q+sP\},(a-sR)_+) &\le J(a)M(s),\\
 N((1-sR)_+,\max\{a,q+s(aP-qR)\}) &\le J(a)M(s),\\
 (-a+sR)_+ &\le J(a)M(s).
\end{aligned}
\end{equation}
\end{definition}

We verify the inequalities in Definition~\ref{def:affine_rows} by computing determinants at the selected anchors.

\begin{lemma}[Inequalities from a determinant perturbation]\label{lem:affine_perturbation_rows}
For every $p\in G$ and every normalized $b$ from
Definition~\ref{def:affine_deficit}, the coordinates $a,q,P,R$ satisfy
the inequalities in Definition~\ref{def:affine_rows}.
\end{lemma}
\begin{proof}
Fix $0<s<1$. Append the single function $sT$ with label $z+A\ell$, where $\ell\in G$.
For $h,\ell\in G$, write $A'_{z+Ah}$ for the new linear part at anchor
$z+Ah$. Define $f_{h,\ell}:=\det(A'_{z+Ah})/\det A$. Set $U:=(h,-s(\ell-h))$ and $V:=(m,b)$, viewed as two-column matrices.
Then
\begin{align}
 f_{h,\ell}
 &=\det(I+h m^\top-s(\ell-h)b^\top)\notag\\*
 &=\det(I_2+V^\top U)\notag\\*
 &=1-\delta_h-s b^\top((1-\delta_h)\ell-(1-\delta_\ell)h),
 \label{eq:affine_rank_two}
\end{align}
where the first step factors $A$ out of the perturbed linear part after the
anchor change and the addition of $sT$, the second applies
Lemma~\ref{lem:affine_sylvester}, and the third expands the $2\times2$
determinant using $m^\top h=-\delta_h$ and $m^\top\ell=-\delta_\ell$.
The ratio of the new denominator in $\Phi$ to its original value is
\[
 (\frac{\int_K e^{D(x)+sT(x)} \d\mu(x)}
                  {\int_K e^{D(x)} \d\mu(x)})^\gamma
 =(\E[e^{sT(X)}])^{\gamma}\le M(s),
\]
where the first step uses the density of $\nu$ and the second is
Lemma~\ref{lem:affine_exponential_moment}. This ratio is independent of the label and the anchors.
By maximality of $\Phi$ and Eq.~\eqref{eq:affine_weight_envelope}, $N$ evaluated at the
positive and negative extrema of the new relative determinants is at most $J(a_+)M(s)$.
This uses symmetry of $N$ if $\det A<0$. Substituting the two labels
$\ell=0,r$ and the three anchors $h=0,p,r$ into Eq.~\eqref{eq:affine_rank_two} gives
\[
\begin{array}{c|ccc}
 & h=0 & h=p & h=r\\ \hline
 \ell=0 & 1 & -q+sP & -a+sR\\
 \ell=r & 1-sR & -q+s(qR-aP) & -a .
\end{array}
\]
The entry $-q+sP$ proves Eq.~\eqref{eq:affine_single_row} because $N$ dominates
each coordinate. For $a\ge0$, the first line of the table bounds the positive
extremum from below by $\max\{1,-q+sP\}$ and the negative extremum from below by
$(a-sR)_+$. The second line bounds these extrema from below by $(1-sR)_+$ and
$\max\{a,q+s(aP-qR)\}$, respectively. By monotonicity of $N$ (Lemma~\ref{lem:affine_monotone}), we obtain the
first two inequalities in Eq.~\eqref{eq:affine_rows}. The single entry at
$\ell=0,h=r$ gives the third.
\end{proof}

\subsection{A common translation}

By Lemma~\ref{lem:affine_perturbation_rows}, the inequalities in Definition~\ref{def:affine_rows} hold for every $p\in G$ and every normalized support direction $b$. We use them to find a single vector $u$ such that $G-u\subseteq C(K-\bar x)$. We take $u=\theta(a)r$, where the coefficient depends only on $a$, so that $u$ is independent of $p$ and $b$.

\begin{definition}[Bands and translations]\label{def:affine_bands}
Set $a_0:=9/20$. For $-1\le a\le a_0$ set $\theta(a):=0$. On $(a_0,1]$ use the
constant $\theta$ on the band containing $a$ in Definition~\ref{def:affine_certificate_data}.
At a shared endpoint take the value from the band on the left.
Each closed band is checked using its own value at both ends.
\end{definition}

For the chosen coefficient $\theta(a)$, we bound $(P-\theta(a)R)$ uniformly over the coordinates $q,P,R$ allowed by Definition~\ref{def:affine_rows}.

\begin{lemma}[Uniform scalar bound]\label{lem:affine_scalar_certificate}
Let $-1\le q\le a\le1$. If $a,q,P,R$ satisfy Definition~\ref{def:affine_rows},
then $P-\theta(a)R<C$.
\end{lemma}
\begin{proof}
If $a\le a_0$, Table~\ref{tab:affine_knots} gives $J(a_+)=1$.
Take $s:=13/20$. Then $M(s)<2$ since $20^{29}<2^{44}7^{29}$. Hence
\[
 P\le(q+M(s))/s<(9/20+2)/(13/20)=49/13<C,
\]
where the first step solves Eq.~\eqref{eq:affine_single_row}, the second uses
$q\le a_0$ and the bound on $M(s)$, and the remaining steps simplify the fraction
and use $49\cdot365<13\cdot1392$.

Suppose $a>a_0$. At $s=1/2$, the table gives $J(1)<5/4$, and
$M(1/2)<8/5$ since $2^{29}5^{44}<8^{44}$. Thus $J(a)M(1/2)<2$.
Eq.~\eqref{eq:affine_rows} gives $1-R/2<2$ and $-a+R/2<2$,
using $N(u,v)\ge u$ for the former. Consequently $-2<R<6$ since $a\le1$.

Fix the band containing $a$ and its value of $\theta$. If $P-\theta R\ge C$,
set $P_0:=C+\theta R$ and $v:=(q+1)/(1+a)$. Then $0\le v\le1$ since $-1\le q\le a$.
Both $-q+sP$ and $q+s(aP-qR)$ are nondecreasing in $P$ since $s>0$ and $a>0$.
By Lemma~\ref{lem:affine_monotone}, Eq.~\eqref{eq:affine_rows} still holds with
$P=P_0$. By Lemma~\ref{lem:affine_finite_coverage} in Appendix~\ref{app:affine_exact},
a row of Eq.~\eqref{eq:affine_rows} fails at some value of $s$. The contradiction proves the claim.
\end{proof}

A bound in every normalized support direction controls membership in $C(K-\bar x)$. The independence of $\theta(a)$ from $p$ and $b$ allows these bounds to apply to the single translate $G-\theta(a)r$.

\begin{lemma}[Common affine placement]\label{lem:affine_common_placement}
Suppose $K,L\subset\R^n$ are convex bodies with $0\in\operatorname{int}K$,
and let $C$ be as in Definition~\ref{def:affine_parameters}. There exist an invertible affine map $F$ with
linear part $A$, points $z\in L$, $\bar x\in\operatorname{int}K$, a vector
$u\in\R^n$, and $44n/29\le H\le73n/29$ such that
\[
 F(K)\subseteq H(L-z),\qquad A^{-1}(L-z)-u\subseteq C(K-\bar x).
\]
\end{lemma}
\begin{proof}
Choose a maximizing list from Lemma~\ref{lem:affine_maximizer} and its
anchor $z$ as in Definition~\ref{def:affine_anchors}. Take $F:=F_z$
and $A:=A_z$. The first inclusion follows from Lemma~\ref{lem:affine_forward_placement}.
The choice of the list and anchors fixes $r$ and $a$ before either a point
$p\in G$ or a support direction $b$ is chosen. The scalar inequalities hold by Lemma~\ref{lem:affine_perturbation_rows}.
Set $u:=\theta(a)r$, which is independent of $p$ and $b$. For each $b$ with
$h_{K-\bar x}(b)=1$, we have $b^\top(p-u)=P-\theta(a)R<C$,
where the equality substitutes $u=\theta(a)r$, $P=b^\top p$ and $R=b^\top r$,
and the inequality is Lemma~\ref{lem:affine_scalar_certificate}.
By Lemma~\ref{lem:affine_support}, the point $(p-u)/C$ belongs to
$K-\bar x$. This proves the second inclusion. The stated range of $H$ follows from
Lemma~\ref{lem:affine_maximizer}.
\end{proof}

\begin{figure}[H]
    \centering
    \definecolor{bmblue}{RGB}{45,112,209}
    \definecolor{bmred}{RGB}{210,55,55}
    \begin{tikzpicture}[x=1.25cm,y=0.85cm,line cap=round,line join=round,>=Latex,font=\small]
        \draw[bmblue,fill=bmblue!5,line width=0.8pt]
            (-2.4,-0.5) -- (-1.6,-1.5) -- (0.7,-1.4) -- (2.2,-0.3)
            -- (1.7,1.3) -- (-0.4,1.7) -- (-2.2,0.8) -- cycle;
        \def\bmtranslatedpolygon{(-1.5,-0.5) -- (-0.4,-0.95) -- (1.25,-0.4) -- (0.8,0.75) -- (-0.8,1.05)}
        \draw[bmred,fill=bmred!8,line width=0.8pt] \bmtranslatedpolygon -- cycle;
        \begin{scope}[shift={(1.65,0.65)}]
            \draw[bmred,dashed,line width=0.8pt] \bmtranslatedpolygon -- cycle;
        \end{scope}
        \draw[->,black,line width=0.8pt]
            (1.505222,0.674063) -- node[pos=0.25,above=3pt,text=black,inner sep=1pt] {$-u$} (-0.144778,0.024063);
        \fill[bmred] (1.505222,0.674063) circle[radius=1.3pt];
        \fill[bmred] (-0.144778,0.024063) circle[radius=1.3pt];
        \node[above,text=bmblue] at (-1.35,1.65) {$C(K-\bar x)$};
        \node[text=bmred] at (-0.45,-0.5) {$G-u$};
        \node[text=bmred] at (2.35,0.85) {$G$};
    \end{tikzpicture}
    \caption{An illustration of the common translation.}
    \label{fig:affine_common_translation}
\end{figure}
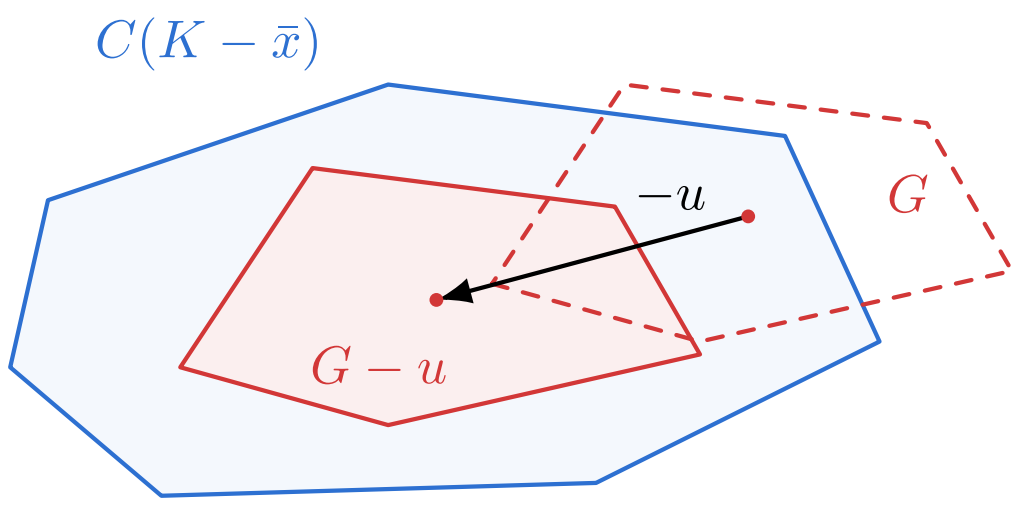

\subsection{Proof of the main theorem}
\label{sec:affine_main_proof}

We prove Theorem~\ref{thm:main} by translating both affine images so that their containments have the form in Definition~\ref{def:affine_distance}.

\begin{proof}
By translation assume $0\in\operatorname{int}K$. Use $C$ from Definition~\ref{def:affine_parameters}.
By Lemma~\ref{lem:affine_common_placement}, there exist an invertible affine map
$F$ with linear part $A$, points $z\in L$, $\bar x\in\operatorname{int}K$,
a vector $u$, and a scalar $H$ with $44n/29\le H\le 73n/29$ such that
\[
 F(K)\subseteq H(L-z),\qquad A^{-1}(L-z)-u\subseteq C(K-\bar x).
\]
Set $B:=F(K)$, $U:=H(L-z)$ and $\lambda:=HC>1$.
Set $w:=HAu-\lambda F(\bar x)$. Then
\[
 B\subseteq U\subseteq HAu+HC A(K-\bar x)=w+\lambda B,
\]
where the first step uses the inner containment in
Lemma~\ref{lem:affine_common_placement}, the second follows by adding $u$
and applying $HA$ to its second inclusion, and the third uses
$A(K-\bar x)=B-F(\bar x)$ and $\lambda=HC$.
With $t:=w/(\lambda-1)$, set $B':=B+t$ and $U':=U+t$. We have
$B'\subseteq U'\subseteq\lambda B'$, where the first inclusion translates
$B\subseteq U$ by $t$, and the second translates $U\subseteq w+\lambda B$
and uses $w+t=\lambda t$. Consequently,
\[
 d_{\mathrm{BM}}(K,L)\le HC\le \frac{73n}{29}\frac{1392}{365}
 =\frac{48}{5}n,
\]
where the first step follows from Definition~\ref{def:affine_distance}
applied to the affine images $B',U'$ with dilation $\lambda=HC$, the second
uses $H\le73n/29$ from Lemma~\ref{lem:affine_common_placement} and
$C=1392/365$, and the third step comes from $1392=48\cdot29$ and $365=5\cdot73$. Taking the supremum over $K,L$ gives the assertion for $D_n$ by Definition~\ref{def:affine_distance}.
\end{proof}

\section*{Acknowledgements and AI Disclosure}
The authors want to thank Brayden Letwin, Victor Reis, and Lichen Zhang for helpful discussions.

The proof was obtained with the assistance of Claude Fable 5.1 and GPT-6 Astra, under the guidance of the authors, who set the overall strategy and decided which directions to pursue. AI tools were also used in writing the paper and in checking the proof, which was further verified by formalizing it in Lean. The authors take full responsibility for the correctness of the paper.

\bibliographystyle{alpha}
\bibliography{ref}

\appendix

\section{Exact verification of the uniform scalar bound}
\label{app:affine_exact}

\subsection{Certificate data}
\label{sec:affine_certificate_data}

The certificate specifies rational knots for $J$, bands in $a$, and one translation coefficient $\theta$ on each band. The constants $\gamma$ and $C$ are fixed in Definition~\ref{def:affine_parameters}. We give the knots below; the bands and their coefficients are generated in Section~\ref{sec:affine_coverage} after defining the interval test.

\begin{definition}[Knots of the determinant function]\label{def:affine_knots}
For each column of each pair of rows in Table~\ref{tab:affine_knots}, set
$x_i:=X_i/200$ and $j_i:=1+Y_i/10^5$. The function $J$ joins the points
$(x_i,j_i)$ by straight segments in increasing order of $x_i$.
Thus $J=1$ on $[0,1/2]$. For $x_i\le z\le x_{i+1}$ its value is
\[
 J(z)=j_i+(z-x_i)(j_{i+1}-j_i)/(x_{i+1}-x_i).
\]
\end{definition}

\begin{table}[ht]
\centering
\begin{tabular}{r@{\quad}rrrrrrrr}
\toprule
$X_i$ & 0&100&101&105&107&109&111&113\\
$Y_i$ & 0&0&492&1253&1389&1897&2315&2658\\
\midrule
$X_i$ & 115&116&119&120&121&123&126&127\\
$Y_i$ & 2957&3015&3015&3239&4099&4555&4555&4628\\
\midrule
$X_i$ & 128&129&131&133&135&136&137&140\\
$Y_i$ & 5451&5885&6357&6357&6672&7435&7842&8612\\
\midrule
$X_i$ & 145&147&149&153&155&161&163&164\\
$Y_i$ & 8612&9285&10771&11659&12082&12082&13380&14075\\
\midrule
$X_i$ & 165&167&168&169&173&175&176&177\\
$Y_i$ & 14285&14653&14918&15119&15119&16293&16552&17001\\
\midrule
$X_i$ & 179&180&181&183&184&185&189&200\\
$Y_i$ & 17486&17490&17980&18464&18608&19085&20191&20197\\
\bottomrule
\end{tabular}
\caption{The 48 knots of $J$. Each column pairs $X_i$ with the $Y_i$ immediately below it.}
\label{tab:affine_knots}
\end{table}

\subsection{The interval test and its soundness}
\label{sec:affine_interval_soundness}

For a fixed band $I$ and translation coefficient $\theta$, we test an interval $E$ of $R$-values at once. We first specify rational bounds and rounding rules so that the comparisons can be made with integers.

\begin{definition}[Rational bounds]\label{def:affine_grid}
Set $S:=10^7$, $[x]_-:=\lfloor Sx\rfloor/S$, $[x]_+:=\lceil Sx\rceil/S$
and $C_*:=[C]_-$, with $C=1392/365$ as in Definition~\ref{def:affine_parameters}.
For $1\le k\le900$ define $s_k:=k/1000$ and
\[
 m_k:=\min\{m\in\mathbb{Z}_{>0}:m^{44}(1000-k)^{29}\ge S^{44}1000^{29}\}.
\]
For a closed interval $I:=[a_-,a_+]\subseteq[9/20,1]$ set
$U_k(I):=[J(a_+)]_+m_k/S$. Then $J(a)M(s_k)\le U_k(I)$ for $a\in I$,
by monotonicity of $J$ and the defining integer inequality for $m_k$.
\end{definition}

\begin{definition}[Interval test]\label{def:affine_interval_test}
Fix $I$ as in Definition~\ref{def:affine_grid}, a translation $\theta\in[0,1]$
and an interval $E:=[r_-,r_+]\subseteq[-2,6]$. For each $k=1,\ldots,900$,
write $s:=s_k$, $U:=U_k(I)$ and set
\[
 u_0:=[(a_--sr_+)_+]_-,\qquad u_1:=[(1-sr_+)_+]_-,
 \qquad u_2:=[(-a_++sr_-)_+]_- .
\]
The inner positive part is $x_+:=\max\{x,0\}$. At a corner
$(a,r)\in\{a_-,a_+\}\times\{r_-,r_+\}$, set $P_*:=C_*+\theta r$ and
$q:=(1+a)v-1$. Define the affine functions of $v$
\[
 F_0(v):=-q+sP_*,\qquad F_1(v):=q+s(aP_*-qr).
\]
For $w=0,1,2$ the following rule gives an interval $V_{w,k}\subseteq[0,1]$,
possibly empty. All endpoints are multiples of $1/S$.
\begin{enumerate}
\item For $w=0,1$, set $d_0:=1$, $d_1:=a_-$. If $N(u_w,d_w)>U$, take
$V_{w,k}:=[0,1]$. Otherwise let $b$ be the smallest nonnegative multiple
of $1/S$ with $N(u_w,b)>U$. It exists because $N(u_w,b)\ge b$.
Intersect $[0,1]$ with the four half-lines $F_w(v)\ge b$, one at each
corner of $I\times E$. Round the lower endpoint up and the upper endpoint
down to multiples of $1/S$. Take $V_{w,k}$ to be the resulting interval,
with the empty interval understood if the endpoints cross or the intersection is empty.
A constant $F_w$ imposes either no restriction or an empty intersection.
\item For $w=2$ take $V_{2,k}:=[0,1]$ if $u_2>U$, and
$V_{2,k}:=\varnothing$ otherwise.
\end{enumerate}
The test accepts $I,\theta,E$ if one of these intervals contains $0$,
another contains $1$, and the two intersect. The same interval may serve twice.
\end{definition}

We next justify using these finite comparisons to exclude a whole parameter box, including its boundary.

\begin{lemma}[Soundness of the interval test]\label{lem:affine_integer_test}
If the interval test in Definition~\ref{def:affine_interval_test} accepts $I,\theta,E$, then for
every $(a,R,v)\in I\times E\times[0,1]$, at least one row of
Eq.~\eqref{eq:affine_rows} fails with $q=(1+a)v-1$ and $P=C+\theta R$.
\end{lemma}
\begin{proof}
Fix $v\in V_{w,k}$. In the first two rows the left-hand side is at least
$N(u_w,d_w)$, by the endpoint bounds on the first coordinate (after using
symmetry in row $0$) and the constant anchor in the other coordinate.
If this already exceeds $U_k(I)$, the row fails. In the remaining case,
the inequalities $F_w(v)\ge b$ hold at the four corners of $I\times E$.
For this fixed $v$, successive linear interpolation in $a$ and $r$
expresses $F_w(v)$ as a convex combination of these corner values.
Thus $F_w(v)\ge b$ for every $(a,r)\in I\times E$. Replacing $C_*$ by $C$ can
only increase $F_w(v)$ because its coefficient of $C_*$ is $s$ or $sa$.
The left-hand side is therefore at least $N(u_w,b)>U_k(I)$ by
Lemma~\ref{lem:affine_monotone}. In row $2$ the lower bound is $u_2$.
In each case, the right-hand side is at most $U_k(I)$
by Definition~\ref{def:affine_grid}. The two intersecting intervals cover $[0,1]$, including every
real $v$ between their endpoints.
\end{proof}

\subsection{Coverage and reproducibility}
\label{sec:affine_coverage}

To apply the interval test throughout $[-2,6]$, we choose a single translation coefficient on each band and cover the $R$-range by accepted intervals. The two recurrences below choose the coefficients and subdivide the bands and $R$-intervals when needed.

\begin{definition}[Bands and their translations]\label{def:affine_certificate_data}
Bisect an interval $[l,h]$ into the closed children $[l,m]$ and $[m,h]$,
where $m:=[(l+h)/2]_-$. Use the following two finite recurrences.
\begin{enumerate}
\item For fixed $I,\theta$, let $T_d(E)$ count intervals accepted in $E$,
with $\bot$ denoting failure. If the interval test accepts $I,\theta,E$, set
$T_d(E):=1$. Otherwise return $\bot$ at depth $d=15$. At smaller depths set
$T_d(E):=T_{d+1}(E_1)+T_{d+1}(E_2)$ for the two children of $E$,
where a sum with $\bot$ is $\bot$.
\item Starting at depth $e=0$, use as initial bands
\[
 I_i:=[9/20+i/2000,9/20+(i+1)/2000],\qquad 0\le i<1100 .
\]
On each band $I$, start with the integer range $[L,H]:=[0,10000]$ and
trial $t:=4000$. For $\theta:=t/10000$, evaluate
$A:=T_0([-2,0])$ and $B:=T_0([0,6])$ using this same $I,\theta$.
If both are integers, accept $I$ with translation $\theta$ and count $A+B$.
If only $A$ fails, replace $H$ by $t-1$. If only $B$ fails, replace
$L$ by $t+1$. In either case, if the range is nonempty, repeat with
$t:=\lfloor(L+H)/2\rfloor$.
If both fail or the range becomes empty, bisect $I$ and apply this rule
afresh to each child at depth $e+1$. If a bisection would be needed at
$e=10$, declare that initial band unsuccessful.
\end{enumerate}
\end{definition}

\begin{figure}[H]
    \centering
    \definecolor{bmblue}{RGB}{45,112,209}
    \definecolor{bmred}{RGB}{210,55,55}
    \begin{tikzpicture}[line cap=round,line join=round,>=Latex,font=\small]
        \fill[bmblue!4] (0,0) rectangle (4.8,2.8);
        \fill[bmred!10] (1.3,0.7) rectangle (3.1,1.75);
        \draw[bmblue!70,line width=0.55pt] (0,0) rectangle (4.8,2.8);
        \draw[bmblue!70,line width=0.55pt] (1.3,0) -- (1.3,2.8) (3.1,0) -- (3.1,2.8);
        \foreach \y in {0.35,0.7,1.225,1.75}
            \draw[bmblue!70,line width=0.55pt] (0,\y) -- (1.3,\y);
        \foreach \y in {0.7,1.75,2.275}
            \draw[bmblue!70,line width=0.55pt] (1.3,\y) -- (3.1,\y);
        \foreach \y in {0.7,1.225,1.75,2.275}
            \draw[bmblue!70,line width=0.55pt] (3.1,\y) -- (4.8,\y);
        \draw[bmred,line width=0.9pt] (1.3,0.7) rectangle (3.1,1.75);
        \node[text=bmred] at (2.2,1.225) {$I\times E$};
        \draw[->,bmblue,line width=0.7pt] (0,0) -- (5.15,0) node[right] {$a$};
        \draw[->,bmblue,line width=0.7pt] (0,0) -- (0,3.18) node[above] {$R$};
        \foreach \y/\t in {0/-2,0.7/0,2.8/6} {
            \draw[bmblue,line width=0.6pt] (-0.06,\y) -- (0.06,\y);
            \node[left,text=bmblue] at (-0.08,\y) {$\t$};
        }
        \node[below,text=bmblue] at (0,-0.08) {$9/20$};
        \node[below,text=bmblue] at (4.8,-0.08) {$1$};
        \node[above,text=bmblue] at (0.65,2.85) {$\theta_1$};
        \node[above,text=bmred] at (2.2,2.85) {$\theta_2$};
        \node[above,text=bmblue] at (3.95,2.85) {$\theta_3$};
        \draw[bmred,line width=0.6pt] (1.3,-0.2) -- (1.3,-0.32) -- (3.1,-0.32) -- (3.1,-0.2);
        \node[below,text=bmred] at (2.2,-0.32) {$I$};        \draw[bmred,line width=0.6pt] (-0.46,0.7) -- (-0.58,0.7) -- (-0.58,1.75) -- (-0.46,1.75);        \node[left=2pt,text=bmred] at (-0.58,1.225) {$E$};
        \draw[->,black!55,line width=0.7pt] (5.2,1.4) -- (6.05,1.4);
        \draw[bmblue!40,dotted,line width=0.55pt] (6.4,0.65) -- (6.4,2.2);
        \draw[bmblue!40,dotted,line width=0.55pt] (10.3,0.65) -- (10.3,1.5);
        \draw[->,bmblue,line width=0.7pt] (6.25,0.65) -- (10.65,0.65) node[right] {$v$};
        \draw[bmblue,line width=0.6pt] (6.4,0.59) -- (6.4,0.71) (10.3,0.59) -- (10.3,0.71);
        \node[below,text=bmblue] at (6.4,0.58) {$0$};
        \node[below,text=bmblue] at (10.3,0.58) {$1$};
        \draw[bmblue,line width=1.2pt] (6.4,2.2) -- (9.15,2.2);
        \fill[bmblue] (6.4,2.2) circle[radius=1.2pt] (9.15,2.2) circle[radius=1.2pt];
        \node[above=3pt,text=bmblue] at (7.775,2.2) {$V_{w,k}$};
        \draw[bmred,line width=1.2pt] (8.15,1.5) -- (10.3,1.5);
        \fill[bmred] (8.15,1.5) circle[radius=1.2pt] (10.3,1.5) circle[radius=1.2pt];
        \node[below=3pt,text=bmred] at (9.225,1.5) {$V_{w',k'}$};
    \end{tikzpicture}
    \caption{A schematic illustration of the certificate coverage.
    Left: rectangles $I\times E$ cover the $(a,R)$-domain. Here $I$ is the horizontal $a$-interval and $E$ is the vertical $R$-interval, so each rectangle has width $|I|$ and height $|E|$.
    Each vertical band $I$ uses one fixed translation coefficient $\theta$,
    and its $R$-range is subdivided into intervals $E$; one rectangle is highlighted in red.
    Right: $V_{w,k}$ and $V_{w',k'}$ consist of $v$-values for which the corresponding
    scalar inequalities fail throughout the highlighted rectangle at $s_k$ and $s_{k'}$,
    respectively (Definition~\ref{def:affine_interval_test}).
    These closed intervals overlap and contain $0$ and $1$, respectively,
    so they cover $[0,1]$ and the interval test accepts the highlighted rectangle.}
    \label{fig:affine_certificate_coverage}
\end{figure}
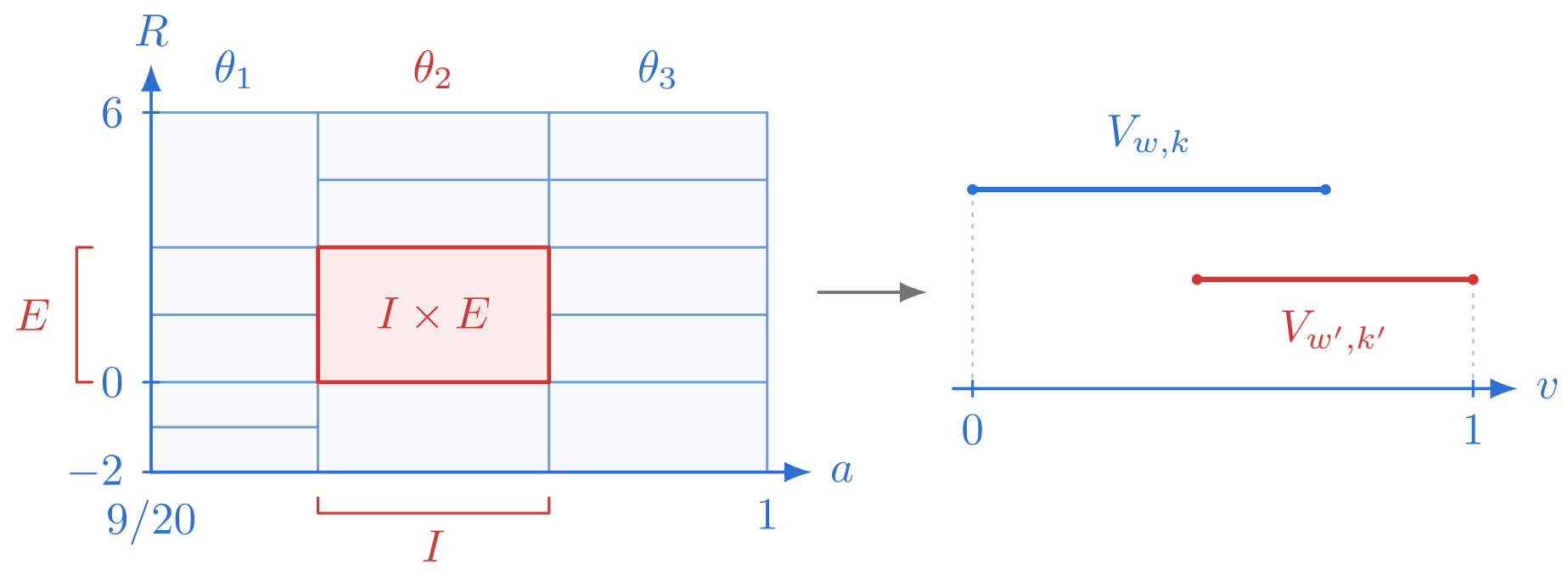

It remains to verify that the recurrences succeed within their depth limits.

\begin{lemma}[Finite coverage]\label{lem:affine_finite_coverage}
All $1100$ initial bands in Definition~\ref{def:affine_certificate_data}
succeed. Their $1507$ accepted bands cover $[9/20,1]$. On each accepted band
one translation suffices in Lemma~\ref{lem:affine_integer_test} for all
$-2\le R\le6$ and $0\le v\le1$.
\end{lemma}
\begin{proof}
The following are the integer evaluations of the two recurrences. The last two
columns give the largest depths $e,d$ on accepted bands and accepted intervals,
respectively. The interval column sums $A+B$ over the accepted bands.
\begin{center}
\begin{tabular}{lrrrrr}
\toprule
Initial index $i$ & Successful roots & Accepted bands & Intervals & $e$ & $d$\\
\midrule
0--199     & 200 & 248 & 27946 & 2 & 15\\
200--399   & 200 & 257 & 30531 & 2 & 15\\
400--599   & 200 & 262 & 19885 & 2 & 15\\
600--799   & 200 & 290 & 22417 & 2 & 15\\
800--999   & 200 & 350 & 20721 & 2 & 15\\
1000--1099 & 100 & 100 &  2112 & 0 & 11\\
\midrule
All        &1100 &1507 &123612 & 2 & 15\\
\bottomrule
\end{tabular}
\end{center}
For clarity, each entry can be evaluated with integer comparisons as follows.
Obtain $m_k$ by bisection between $0$ and $8S$ in
Definition~\ref{def:affine_grid}. The upper endpoint suffices since
$8^{44}>10^{29}$ and $1000-k\ge100$. For $b$ in the interval test, bisect
the integer indices from $0$ to $\lfloor SU\rfloor+1$, using monotonicity
of $N$. To perform a comparison $N(u,y)>U$, set
$B':=\max\{u,y\}$, $b':=\min(u,y)$. If $B'=0$ the comparison fails.
Otherwise choose consecutive knots with $x_iB'\le b'\le x_{i+1}B'$ and use
\[
 B'j_i(x_{i+1}-x_i)+(b'-x_iB')(j_{i+1}-j_i)>U(x_{i+1}-x_i).
\]
This follows by substituting the interpolation formula into $B'J(b'/B')$
and clearing the positive knot spacing. To test intersection, take the largest
right endpoint among the intervals containing $0$, and the smallest left
endpoint among those containing $1$. Both must exist and the former must
be at least the latter. This determines each leaf count before the
additions in Definition~\ref{def:affine_certificate_data}.

Each bisection covers its parent, so a accepted band with integer $A,B$
covers both initial $R$-intervals by induction up their subdivision trees.
The same induction on the band subdivisions gives coverage of every $I_i$.
On each accepted band the translation was fixed before testing either
$R$-interval or choosing a row for $v$. Lemma~\ref{lem:affine_integer_test}
therefore applies with that translation throughout the band.
\end{proof}

This completes the finite verification used in the $a>9/20$ case of Lemma~\ref{lem:affine_scalar_certificate}. We give an explicit input on which the arithmetic of the interval test can be checked directly.

\begin{lemma}[A worked interval]
On the first initial band $I_0=[9/20,901/2000]$, take $\theta:=1999/10000$.
For $-1\le R\le-1/2$ the rows $w=0,1$ in Definition~\ref{def:affine_interval_test}
suffice below and above $v=112/125$, respectively.
\end{lemma}
\begin{proof}
Table~\ref{tab:affine_knots} gives $J(901/2000)=1$.
The following quantities use $k=560$ for $0\le v\le112/125$ and $k=313$
for $112/125\le v\le1$. In each row $f$ is the minimum of $F_w(v)$ over
the eight endpoint triples of $a,R,v$.
At each corner $(a,R)$, affinity in $v$ gives $F_w(v)\ge f$
throughout the indicated interval.
\[
\begin{array}{c|r|c|c|c|r}
 w&k&u_w&f&N(u_w,[f]_-)&m_k\\ \hline
 0&560&73/100&107754951/62500000&2155099/1250000&17178970\\
 1&313&2313/2000&86910629031/10^{11}&6422013489/(5\cdot10^9)&12807382
\end{array}
\]
The entries in the fourth column follow by substituting the corner
coordinates into $F_w$ with $C_*=\lfloor 10^7 C\rfloor/10^7$. The next
column uses linear interpolation at the smaller coordinate divided by the
larger one. In each row the fifth entry is strictly greater than $m_k/S$.
Thus $[f]_-$ is at least the threshold $b$ whenever that threshold is
needed. The cut $112/125$ lies on the grid, so the two closed intervals
meet there and the test accepts.
\end{proof}

\end{document}